\documentclass[12pt,a4paper,reqno]{amsart}  

\usepackage{mathtools}
\usepackage{url}
 
\usepackage[headings]{fullpage}
\usepackage{comment} 
\usepackage{datetime}  
\usepackage[english]{babel} 

\usepackage[full]{textcomp}
\usepackage[osf]{newtxtext} 
\usepackage{cabin} 
\usepackage[bigdelims,vvarbb]{newtxmath} 
\usepackage[cal=boondoxo]{mathalfa} 
\usepackage{zlmtt}

\usepackage{microtype} 

\newcommand{\eps}{\varepsilon}

\newcommand{\dx}{\mathrm{d}} 
\newcommand{\E}{\mathcal{E}} 
\newcommand{\I}{\mathcal{I}} 
\newcommand{\J}{\mathcal{J}} 

\newcommand{\N}{\mathbb{N}}

\newcommand{\Etilde}{\widetilde{E}}

\newcommand{\Stilde}{\widetilde{S}}
 
\newcommand{\Vtilde}{\widetilde{V}}

\newcommand{\Odi}[1]{\Odip{}{#1}}
\newcommand{\Odig}[1]{\mathcal{O}\Bigl(#1\Bigr)}
\newcommand{\Odim}[1]{\mathcal{O}\bigl(#1\bigr)}   
\newcommand{\Odip}[2]{\mathcal{O}_{#1}\left(#2\right)}
\newcommand{\Odipg}[2]{\mathcal{O}_{#1}\Bigl(#2\Bigr)}  
\newcommand{\Odipm}[2]{\mathcal{O}_{#1} (#2)}  
\newcommand{\odip}[2]{{o}_{#1}\left(#2\right)}
\newcommand{\odi}[1]{\odip{}{#1}}

\allowdisplaybreaks

\renewcommand{\qedsymbol}{$\square$}
\newenvironment{Proof}[1][Proof]{\par\noindent\textbf{#1.}~}
{\hfill\qedsymbol\smallskip\par}

\newtheoremstyle{sltheorems}
{10pt}
{6pt}
{\slshape}
{}
{\bfseries}
{.}
{.5em}
{\thmname{#1}\thmnumber{ #2}\thmnote{ (#3)}}

\theoremstyle{sltheorems} 
\newtheorem{Theorem}{Theorem}
\newtheorem{Lemma}{Lemma} 

\newtheoremstyle{remark}
{10pt}
{6pt}
{\rm} 
{}
{\bfseries}
{.}
{.5em}
{\thmname{#1}\thmnumber{ #2}\thmnote{ (#3)}}
 \theoremstyle{remark}

\makeatletter

\def\env@Biggcases{%
  \let\@ifnextchar\new@ifnextchar
  \Biggl\lbrace
  \def\arraystretch{1.2}%
  \array{@{}l@{\quad}l@{}}%
}
\makeatother

\allowdisplaybreaks

\begin{document} 

\title[Sums of four prime cubes]{Sums of four prime cubes in short intervals} 
\date{}
\author{Alessandro Languasco \lowercase{and} Alessandro Zaccagnini}
%

\begin{abstract} 
We prove that a suitable
asymptotic formula for the average number  of representations
of  integers $n=p_{1}^{3}+p_{2}^{3}+p_{3}^{3}+p_{4}^{3}$,
where  $p_1,p_2,p_3,p_4$ are prime numbers,
holds in intervals shorter than the the ones previously known.
\end{abstract} 
\maketitle
 
\section{Introduction}
Let $N$ be a sufficiently large integer and $1\le H \le N$ an  integer.
Let
\begin{equation}
\label{r-def} 
\sum_{n=N+1}^{N+H} 
r (n),
\quad
\textrm{where}
\quad 
r (n) =  
\sum_{n=p_{1}^{3}+p_{2}^{3}+p_{3}^{3}+p_{4}^{3}}
\!\!\!\!\!\!\!\!
\log p_{1} \log p_{2} \log p_{3}\log p_{4},
\end{equation} 
be a suitable short interval average of the number of representation
of an integer as a sum of four prime cubes.
The problem of representing integers as
sum of prime cubes is quite an old one; we recall that Hua \cite{Hua1938}-\cite{Hua1965}
stated that almost all positive integers  satisfying some necessary congruence
conditions are the sum of five prime cubes  and that Daveport
\cite{Davenport1939} proved that almost all positive integers are the sum of four 
positive cubes. More recent results on the positive proportions
of integers that are the sum of four prime cubes were obtained
by Roth \cite{Roth1951}, Ren \cite{Ren2003} and Liu \cite{Liu2012}.
In fact, see Br\"udern \cite{Brudern1995}, it is conjectured 
that all sufficiently large integers satisfying some necessary congruence
conditions are the sum of four prime cubes.
Here we prove that
\begin{Theorem}
\label{thm-uncond}
 Let $N\ge 2$, $1\le H \le N$ be integers.
Then, for every $\eps>0$, there exists $C=C(\eps)>0$ such that
\begin{align*} 
\sum_{n=N+1}^{N+H} &
r (n) 
=
\Gamma \Bigl(\frac{4}{3}\Bigr)^3
H  N^{1/3} 
+
\Odig{ H  N^{1/3}   \exp \Bigl( -C \Bigl( \frac{ \log N}{\log  \log N} \Bigr)^{1/3} \Bigr) }
 \quad \text{as} \ N \to \infty,
\end{align*}
uniformly for    $N^{13/18 +\eps}\le H \le N^{1-\eps}$, where $\Gamma$ is Euler's
function.
\end{Theorem}

This should be compared with a recent result about the positive proportion
of such integers in short intervals by Liu-Zhao \cite{LiuZ2017}
which holds   for $H=N^{17/18}$.
As an immediate consequence of Theorem \ref{thm-uncond}
we can say that, for $N$ sufficiently large, every interval of
size larger than $N^{13/18 +\eps}$ contains the expected amount of
integers which are a sum of four prime cubes.
We remark that this level is essentially optimal given the known density
estimates for the non trivial zeroes of the Riemann zeta function.
Assuming the Riemann Hypothesis (RH)  holds we can further 
improve the size of $H$.
\begin{Theorem}
\label{thm-RH}
 Let $\eps>0$, $N\ge 2$, $1\le H \le N$ be integers and
 assume the Riemann Hypothesis (RH) holds.
Then there exists a   constant $B>3/2$ such that
\begin{align*} 
\sum_{n=N+1}^{N+H} &
r (n) 
=
\Gamma \Bigl(\frac{4}{3}\Bigr)^3
H  N^{1/3} 
+
\Odig{ \frac{H^2}{N^{2/3}}
+ H^{3/4}N^{5/12+\eps}
+ H^{1/2} N^{2/3}(\log N)^B
+N(\log N)^3}
\end{align*}
as $N \to \infty$,
uniformly for $\infty(N^{2/3}L^{2B}) \le H \le \odi{N}$,
where $f=\infty(g)$ means $g=\odi{f}$ and $\Gamma$ is Euler's
function.
\end{Theorem}

As an immediate consequence of Theorem \ref{thm-RH}
we can say that, for $N$ sufficiently large, every interval of
size larger than $N^{2/3+\eps}$ contains  the expected amount of
integers which are a sum of four prime cubes. 
We remark that this level is essentially optimal given the spacing
of the cubic sequence.
%

In both the proofs of Theorems \ref{thm-uncond}-\ref{thm-RH} we will use the original
Hardy-Littlewood generating functions to exploit
the easier main term treatment they allow (comparing
with the one which would follow
using  Lemmas 2.3 and 2.9 of Vaughan \cite{Vaughan1997}).

\section{Setting}
Let   $e(\alpha) = e^{2\pi i\alpha}$,
$\alpha \in [-1/2,1/2]$, $L=\log N$, $z= 1/N-2\pi i\alpha$,
\begin{equation*} 
\Stilde_\ell(\alpha) = \sum_{n=1}^{\infty} \Lambda(n) e^{-n^{\ell}/N} e(n^{\ell}\alpha)   \quad
\textrm{and}
\quad
\Vtilde_\ell(\alpha) = \sum_{p=2}^{\infty} \log p \, e^{-p^{\ell}/N} e(p^{\ell}\alpha).
\end{equation*} 
We remark that 
\begin{equation}
\label{z-estim}
\vert z\vert^{-1} \ll \min \bigl(N, \vert \alpha \vert^{-1}\bigr).
\end{equation}
We further set
\begin{equation}
\notag 
   U(\alpha,H)
  = 
  \sum_{m=1}^H e(m \alpha)  
\end{equation}
and, moreover, we also have the usual 
numerically explicit inequality
\begin{equation}
\label{UH-estim}
\vert U(\alpha,H) \vert
\le
\min \bigl(H; \vert \alpha \vert^{-1}\bigr),
\end{equation}
see, \emph{e.g.}, on page 39 of Montgomery \cite{Montgomery1994}. %
We list now the needed preliminary results.
%

 \begin{Lemma}[Lemma 3 of  \cite{LanguascoZ2016b}]
\label{tilde-trivial-lemma}
Let $\ell\ge 1$ be an integer. Then
\(
\vert \Stilde_{\ell}(\alpha)- \Vtilde_{\ell}(\alpha) \vert 
\ll_{\ell}
 N^{1/(2\ell)}  .
\)
\end{Lemma}

\begin{Lemma}
\label{Linnik-lemma}
Let $\ell \ge 1$ be an integer, $N \ge 2$  and $\alpha\in [-1/2,1/2]$.
Then
\begin{equation*}
\Stilde_{\ell}(\alpha)  
= 
\frac{\Gamma(1/\ell)}{\ell z^{1/\ell}}
- 
\frac{1}{\ell}\sum_{\rho}z^{-\rho/\ell}\Gamma (\rho/\ell ) 
+
\Odip{\ell}{1},
\end{equation*}
where $\rho=\beta+i\gamma$ runs over
the non-trivial zeros of $\zeta(s)$.
\end{Lemma}
\begin{Proof}
It follows the line of  Lemma 2 of  \cite{LanguascoZ2016a}; we just 
correct an oversight in its  proof. In eq. (5) on page 48 of 
\cite{LanguascoZ2016a} the term 
\(
-  
  \sum_{m=1}^{\ell \sqrt{3}/4} \Gamma (- 2m/\ell ) z^{2m/\ell}
\)
is missing. Its estimate is trivially $\ll_{\ell} \vert z \vert^{\sqrt{3}/2} \ll_{\ell} 1$.
Hence such an oversight does not affect the final result of  
Lemma 2 of  \cite{LanguascoZ2016a}.
\end{Proof}

\begin{Lemma} [Lemma 4 of \cite{LanguascoZ2016a}]
 \label{Laplace-formula}
Let $N$ be a positive integer and 
$\mu > 0$.
Then
\[
  \int_{-1 / 2}^{1 / 2} z^{-\mu} e(-n \alpha) \, \dx \alpha
  =
  e^{- n / N} \frac{n^{\mu - 1}}{\Gamma(\mu)}
  +
  \Odipg{\mu}{\frac{1}{n}},
\]
uniformly for $n \ge 1$.
\end{Lemma}

\begin{Lemma}
 \label{LP-Lemma-gen} 
Let $\eps$ be an arbitrarily small
positive constant,  $\ell \ge 1$ be an integer, $N$ be a
sufficiently large integer and $L= \log N$. Then there exists a positive constant 
$c_1 = c_{1}(\eps)$, which does not depend on $\ell$, such that 
\[
\int_{-\xi}^{\xi} \,
\Bigl\vert
\Stilde_\ell(\alpha) - \frac{\Gamma(1/\ell)}{\ell z^{1/\ell}}
\Bigr\vert^{2}
\dx \alpha 
\ll_{\ell}
 N^{2/\ell-1} \exp \Big( - c_{1}  \Big( \frac{L}{\log L} \Big)^{1/3} \Big) 
\]
uniformly for $ 0\le \xi < N^{-1 +5/(6\ell) - \eps}$.
Assuming RH we get 
\[
\int_{-\xi}^{\xi} \,
\Bigl\vert
\Stilde_\ell(\alpha) - \frac{\Gamma(1/\ell)}{\ell z^{1/\ell}}
\Bigr\vert^{2}
\dx \alpha 
\ll_{\ell}
N^{1/\ell}\xi L^{2}
\]
uniformly  for  $0 \le \xi \le 1/2$.
\end{Lemma} 
\begin{Proof}
It follows the line of Lemma 3 of \cite{LanguascoZ2016a} and 
Lemma 1 of \cite{LanguascoZ2016b}; we just 
correct an oversight in their  proofs which is based on Lemma \ref{Linnik-lemma} above. Both eq. (8) on page 49 of 
\cite{LanguascoZ2016a}  and eq. (6) on page 423
of \cite{LanguascoZ2016b} should read as 
\[
\int_{1/N}^{\xi}
\Big \vert\sum_{\rho\colon \gamma > 0}z^{-\rho/\ell}\Gamma (\rho/\ell ) \Big \vert^2 
\dx \alpha 
\le
\sum_{k=1}^K
\int_\eta^{2\eta} \Big \vert\sum_{\rho\colon \gamma > 0}z^{-\rho/\ell}\Gamma (\rho/\ell ) \Big \vert^2 \dx \alpha,
\]
where $\eta=\eta_k= \xi/2^k$, $1/N\le \eta \le \xi/2$  and $K$ is a suitable integer satisfying $K=\Odi{L}$. 
The remaining part of the proofs are left untouched. 
Hence such  oversights do not affect the final result of  
Lemma 3 of \cite{LanguascoZ2016a} and 
Lemma 1 of \cite{LanguascoZ2016b}.
\end{Proof}

In the unconditional case a crucial role is played by 
\begin{Lemma}[Hua]
\label{Hua-lemma-series}
Let $N$ be sufficiently large, $\ell,k$ integers, $\ell \ge 1$, $1\le k \le \ell$.
There exists a suitable positive constant $A=A(\ell,k)$ such that
\[ 
\int_{-1/2}^{1/2} 
\vert \Stilde_{\ell}(\alpha)\vert ^{2^k} \ \dx\alpha 
\ll_{k,\ell} 
 N^{(2^k-k)/\ell} L^A
 \quad
 \text{and}
 \quad
 \int_{-1/2}^{1/2} 
\vert \Vtilde_{\ell}(\alpha)\vert^{2^k} \ \dx\alpha 
\ll_{k,\ell} 
 N^{(2^k-k)/\ell} L^A.
 \]
\end{Lemma} 
\begin{Proof}
We  just prove the first part since the second one follows immediately by remarking that
 the   primes are supported on a thinner set than   the prime powers.
Let  $P=(2NL/\ell)^{1/\ell}$. A direct estimate
gives   $\Stilde_{\ell}(\alpha)= \sum_{n\le P} \Lambda(n) e^{-n^\ell/N} e(n^\ell\alpha) + \Odipm{\ell}{L^{1/\ell}}$.
Recalling that the Prime Number Theorem implies
$S_{\ell}(\alpha;t) := \sum_{n\le t} \Lambda(n)   e(n^{\ell}\alpha)  \ll t$, a
partial integration argument gives 
\[
\sum_{n\le P} \Lambda(n) e^{-n^\ell/N} e(n^\ell\alpha) 
= 
-\frac{\ell}{N} \int_1^P t^{\ell-1} e^{-t^\ell/N} S_{\ell}(\alpha;t) \  \dx t
+ \Odipm{\ell}{L^{1/\ell}}.
\]
Using the inequality $(\vert a\vert + \vert b \vert)^{2^k} \ll_k \vert a\vert^{2^k} + \vert b \vert^{2^k}$, H\"older's inequality
and interchanging the integrals, we get that
\begin{align*}
 \int_{-1/2}^{1/2} 
&\vert \Stilde_{\ell}(\alpha)\vert ^{2^k} \ \dx\alpha 
\ll_{k,\ell} 
 \int_{-1/2}^{1/2} 
\Bigl\vert \frac{1}{N} \int_1^P t^{\ell-1} e^{-t^\ell/N} S_{\ell}(\alpha;t)\  \dx t \Bigr\vert ^{2^k} \ \dx\alpha 
+  L^{2^k/\ell}
\\
&\ll_{k,\ell} 
\frac{1}{N^{2^k}}
\Bigl( \int_1^P t^{\ell-1} e^{-t^\ell/N} \  \dx t \Bigr)^{2^k-1}
\Bigl( \int_1^P t^{\ell-1} e^{-t^\ell/N} \int_{-1/2}^{1/2} \vert S_{\ell}(\alpha;t)\vert ^{2^k} \dx\alpha  \  \dx t \Bigr)
+  L^{2^k/\ell}.
\end{align*} 
Theorem 4 of Hua \cite{Hua1965} implies, remarking that
the von Mangoldt function is supported on a thinner set than the integers and
inserts a logarithmic weight whose total contribution can be inserted in the
power of $L$,  that there exists a positive constant $B=B(\ell,k)$ such that
$\int_{-1/2}^{1/2} \vert S_{\ell}(\alpha;t)\vert^{2^k}\ \dx\alpha \ll_{k,\ell} t^{2^k-k} (\log t)^B$.
Using such an estimate 
and remarking that $ \int_1^P t^{\ell-1} e^{-t^\ell/N} \  \dx t \ll_\ell N$,
we obtain  that
\begin{align*}
 \int_{-1/2}^{1/2} 
\vert \Stilde_{\ell}(\alpha)\vert ^{2^k} \ \dx\alpha 
&\ \ll_{k,\ell} 
\frac{1}{N} 
 \int_1^P t^{\ell-1+2^k-k} e^{-t^\ell/N} (\log t)^B \  \dx t  
+  L^{2^k/\ell}
 \ll_{k,\ell}   N^{(2^k-k)/\ell} L^{B+(2^k-k)/\ell} 
\end{align*}
by a direct computation.
This proves the first part of the lemma. \end{Proof}

In fact the argument used in the proof of Lemma \ref{Hua-lemma-series}
can  be used to derive other estimates on $\Stilde_{\ell}(\alpha)$ from the 
ones on  $S_{\ell}(\alpha;t)$. Another instance
of this fact is the following  lemma   about the truncated fourth-mean average of 
$ \Stilde_{\ell}(\alpha)$ which is based on a result by Robert-Sargos \cite{RobertS2006}.

\begin{Lemma}
\label{S-ell_quarta_tau}
Let $N\in \N$, $\eps>0$,   $\ell>1$ and $\tau>0$. Then we have
\[
\int_{-\tau}^{\tau}\vert \Stilde_{\ell}(\alpha)\vert ^4\, \dx \alpha
\ll
\bigl(\tau N^{2/\ell}+N^{4/\ell-1}\bigr)N^{\eps}
 \quad
 \text{and}
 \quad
\int_{-\tau}^{\tau}\vert \Vtilde_{\ell}(\alpha)\vert ^4\, \dx \alpha
\ll
\bigl(\tau N^{2/\ell}+N^{4/\ell-1}\bigr)N^{\eps}.
\]
\end{Lemma}
\begin{Proof}
We can argue   as in the proof of Lemma \ref{Hua-lemma-series}
using  
Lemma 4 of \cite{GambiniLZ2018} on 
$S_{\ell}(\alpha;t) = \sum_{n\le t} \Lambda(n) e(n^\ell\alpha)$
instead of  Theorem 4 of Hua \cite{Hua1965}.
\end{Proof}

The last lemma is a consequence of Lemma \ref{S-ell_quarta_tau}.
\begin{Lemma}
\label{S-ell_quarta_coda}
Let $N\in \N$, $\eps>0$,   $\ell>2$,   $c\ge 1$ and $N^{-c}\le \omega \le N^{2/\ell-1}$. 
Then we have
\[
\Bigl(\int_{-1/2}^{-\omega}
+\int_{\omega}^{1/2}
\Bigr)
\vert \Stilde_{\ell}(\alpha)\vert ^4 \frac{\dx \alpha}{\vert \alpha\vert}
\ll
\frac{ N^{4/\ell-1+\eps}}{\omega}
 \quad
 \text{and}
 \quad
\Bigl(\int_{-1/2}^{-\omega}
+\int_{\omega}^{1/2}
\Bigr)
\vert \Vtilde_{\ell}(\alpha)\vert ^4 \frac{\dx \alpha}{\vert \alpha\vert}
\ll
\frac{ N^{4/\ell-1+\eps}}{\omega}.
\]
\end{Lemma}
\begin{Proof}
By partial integration and Lemma \ref{S-ell_quarta_tau} we get that
\begin{align*}
\int_{\omega}^{1/2}
\vert \Stilde_{\ell}(\alpha)\vert ^4 \frac{\dx \alpha}{ \alpha}
&\ll
\frac{1}{\omega}
\int_{-\omega}^{\omega} 
\vert \Stilde_{\ell}(\alpha) \vert^{4}\ \dx \alpha 
+
\int_{-1/2}^{1/2} 
\vert \Stilde_{\ell}(\alpha) \vert^{4}\ \dx \alpha 
+
\int_{\omega}^{1/2} 
\Bigl(
\int_{-\xi}^{\xi} 
\vert \Stilde_{\ell}(\alpha) \vert^{4}\ \dx \alpha 
\Bigr)
\frac{\dx \xi}{\xi^2}
\\&
\ll
\frac{1}{\omega}\bigl(\omega N^{2/\ell}+N^{4/\ell-1}\bigr)N^{\eps}
+N^{2/\ell+\eps}
+
N^{\eps}
\int_{\omega}^{1/2} 
\frac{\xi N^{2/\ell}+N^{4/\ell-1}}{\xi^2}\dx \xi
\\&
\ll 
N^{2/\ell+\eps} \vert \log (2\omega) \vert
+
\frac{ N^{4/\ell-1+\eps}}{\omega}
\ll
\frac{ N^{4/\ell-1+\eps}}{\omega}
\end{align*}
since $N^{-c}\le \omega \le N^{2/\ell-1}$.
A similar computation proves the result in $[-1/2,-\omega]$ too.
The estimate on $\Vtilde_{\ell}(\alpha)$ can be obtained analogously.
\end{Proof}

 \section{The unconditional case}
 \label{unconditional}
Let $H>2B$, where
\begin{equation}
\label{B-def} 
B= N^{2\eps}.
\end{equation} 
Letting $I(B,H):=[-1/2,-B/H]\cup  [B/H, 1/2]$,
and  recalling \eqref{r-def},
we have  
\begin{align}
 \notag
\sum_{n=N+1}^{N+H} e^{-n/N}  r(n)
& =  
\int_{-1/2}^{1/2}\Vtilde_{3}(\alpha)^4    U(-\alpha,H)e(-N\alpha) \, \dx \alpha  
 \\
 \notag
 & 
 =   
\int_{-B/H}^{B/H}    \Stilde_{3}(\alpha)^4 
 U(-\alpha,H)e(-N\alpha) \, \dx \alpha 
 +
\int\limits_{I(B,H)}    \Stilde_{3}(\alpha)^4  
 U(-\alpha,H)e(-N\alpha) \, \dx \alpha 
 \\
 \label{main-dissection-series}
 & \hskip1cm
 +
 \int_{-1/2}^{1/2} \bigl( \Vtilde_{3}(\alpha)^4  -\Stilde_{3}(\alpha)^4 \bigr)  
 U(-\alpha,H)e(-N\alpha) \, \dx \alpha
 =
 I_1+I_2+I_3,
 \end{align} 
 say.
 Now we evaluate these terms.
 
 \subsection{Estimation of $I_2$}
%
 Using \eqref{UH-estim}
and  Lemma \ref{S-ell_quarta_coda} with $\omega = B/H$ and $\ell=3$, we obtain 
\begin{equation}
\label{I2-estim}
I_2
\ll 
\int_{B/H}^{1/2} 
\vert \Stilde_{3}(\alpha) \vert^{4} \frac{\dx \alpha}{\alpha}
 \ll  
\frac{H }{B }  N^{1/3+\eps},
\end{equation}
provided that $H\gg N^{1/3}B$.
 
\subsection{Estimation of $I_3$}
\label{I3-estim}
Clearly
\begin{align*}
\vert \Vtilde_{3}(\alpha)^4  -\Stilde_{3}(\alpha)^4 \vert 
&=
\vert\Vtilde_{3}(\alpha)  -\Stilde_{3}(\alpha) \vert
 \vert \Vtilde_{3}(\alpha)^3+\Vtilde_{3}(\alpha)^2\Stilde_{3}(\alpha)  +\Vtilde_{3}(\alpha) \Stilde_{3}(\alpha)^2+\Stilde_{3}(\alpha)^3\vert
\\
& \ll
 \vert\Vtilde_{3}(\alpha)  -\Stilde_{3}(\alpha) \vert
  (  \vert \Vtilde_{3}(\alpha)  \vert + \vert \Stilde_{3}(\alpha) \vert ) ^3
\\
&  \ll
  \vert\Vtilde_{3}(\alpha)  -\Stilde_{3}(\alpha) \vert
 \max (  \vert \Vtilde_{3}(\alpha)\vert ^3 ;  \vert \Stilde_{3}(\alpha)\vert^3  ).
\end{align*}

Hence by  Lemma \ref{tilde-trivial-lemma} we have
\begin{equation}
 \label{I3-estim-1}
I_3
\ll
N^{1/6}
\int_{-1/2}^{1/2} 
\bigl(
\vert \Vtilde_{3}(\alpha)\vert^3 
+
\vert \Stilde_{3}(\alpha)\vert^3  
\bigr)
\vert U(-\alpha,H)
\vert \, \dx \alpha 
= N^{1/6}(K_1+K_2),
\end{equation}
say.
%
  %
By \eqref{UH-estim} we get
\begin{equation}
\label{K2-split}
K_2
\ll 
H
\int_{-1/H}^{1/H} 
 \vert \Stilde_{3}(\alpha)\vert^3 \, \dx \alpha 
 + 
 \Bigl(
\int_{-1/2}^{-1/H} + 
\int_{1/H}^{1/2} 
\Bigr)
 \vert \Stilde_{3}(\alpha)\vert^3  \frac{\dx \alpha}{\vert \alpha\vert}
 = K_{2,1}+K_{2,2}, 
 \end{equation}
 say.
Using the H\"older inequality and Lemma  \ref{S-ell_quarta_tau} with $\tau = 1/H$ and $\ell=3$
we get 
\begin{equation}
\label{K21-estim}
K_{2,1}
\ll
H^{3/4}
 \Bigl(
\int_{-1/H}^{1/H} 
 \vert \Stilde_{3}(\alpha)\vert^4 \, \dx \alpha 
 \Bigr)^{3/4}
\ll
H^{3/4}N^{1/4+\eps},
\end{equation}
provided that $H\gg N^{1/3}$.
Using the H\"older inequality
and Lemma \ref{S-ell_quarta_coda} with $\omega = 1/H$ and $\ell=3$
we get 
\begin{equation}
 \label{K22-estim}
K_{2,2}
\ll
 \Bigl(
\int_{1/H}^{1/2} 
 \vert \Stilde_{3}(\alpha)\vert^4 \frac{\dx \alpha}{\alpha}
 \Bigr)^{3/4}
  \Bigl(
\int_{1/H}^{1/2} 
\frac{\dx \alpha}{\alpha}
 \Bigr)^{1/4} 
 \ll
  \bigl(
H N^{1/3+\eps}
 \bigr)^{3/4}
 L^{1/4}
 \ll
H^{3/4}N^{1/4+\eps},
\end{equation} 
provided that $H\gg N^{1/3}$.
Combining \eqref{K2-split}-\eqref{K22-estim} we obtain
\begin{equation}
\label{K2-estim}
K_2
\ll
H^{3/4}N^{1/4+\eps},
\end{equation}
provided that $H\gg N^{1/3}$.
An analogous computation gives
  %
\begin{equation}
\label{K1-estim}
K_1
\ll
H^{3/4}N^{1/4+\eps},
\end{equation}
and, by \eqref{I3-estim-1} and \eqref{K2-estim}-\eqref{K1-estim}, we can finally write
\begin{equation}
 \label{I3-estim-series}
I_3
\ll  
H^{3/4}N^{5/12+\eps},
\end{equation}
provided that $H\gg N^{1/3}$.

 \subsection{Evaluation of $I_1$}
 
 Since $\Gamma(4/3)= (1/3)\Gamma(1/3)$,
 we have that
 \begin{align}
\notag
I_1 &=  
\int_{-B/H}^{B/H}     \frac{\Gamma(4/3)^4}{z^{4/3}} 
 U(-\alpha,H)e(-N\alpha) \, \dx \alpha 
 +
\int_{-B/H}^{B/H}    \Bigl(\Stilde_{3}(\alpha)^4   - \frac{\Gamma(4/3)^4}{z^{4/3}} \Bigr)
 U(-\alpha,H)e(-N\alpha) \, \dx \alpha 
 \\&
  \label{I1-split}
 = J_1+J_2,
 \end{align}
say.
By \eqref{z-estim}-\eqref{UH-estim} and Lemma \ref{Laplace-formula},
a direct calculation  gives
\begin{align}
\notag
J_1
&
=
\Gamma\Bigl(\frac{4}{3}\Bigr)^3\
\sum_{n=N+1}^{N+H} e^{-n/N} n^{1/3}
+\Odig{\frac{H}{N}}  +
\Odig{\int_{B/H}^{1/2} \frac{\dx \alpha}{\alpha^{7/3}}}
\\
\notag &
=
\frac{\Gamma(4/3)^3}{e}\
\sum_{n=N+1}^{N+H} n^{1/3}
+\Odig{\frac{H}{N}+\frac{H^2}{N^{2/3}} + \frac{H^{4/3}}{B^{4/3}}}
\\
\label{J1-eval-series}
&
=
\Gamma\Bigl(\frac{4}{3}\Bigr)^3
\frac{HN^{1/3}}{e}
+\Odig{ \frac{H^{4/3}}{B^{4/3}}+N^{1/3}}.
 \end{align}
From now on, we denote 
\(
\Etilde_{3}(\alpha) : =\Stilde_3(\alpha) - \frac{\Gamma(4/3)}{z^{1/3}}.
\)
By  $f^2 - g^2 = 2g (f-g) + (f-g)^2$, \eqref{z-estim}  and $\Stilde_{3}(\alpha) \ll N^{1/3}$
we get 
\begin{align}
\notag
\Stilde_{3}(\alpha)^4  - \frac{\Gamma(4/3)^4}{z^{4/3}} 
&=
\Bigl( 
\Stilde_{3}(\alpha) ^2  +  
\frac{\Gamma(4/3)^2}{z^{2/3}}
\Bigr)
\Bigl(
\Stilde_{3}(\alpha) ^2  - 
\frac{\Gamma(4/3)^2}{z^{2/3}}
\Bigr)
\\ 
\notag
&=
\Bigl( 
\Stilde_{3}(\alpha) ^2  +  
\frac{\Gamma(4/3)^2}{z^{2/3}}
\Bigr)
\Bigl( 
2\frac{\Gamma(4/3)}{z^{1/3}}   \Etilde_{3}(\alpha)  
+
  \Etilde_{3}(\alpha)  ^2
\Bigr)
\\ 
\label{S-tilde-uncond-approx} 
& \ll
\vert \Stilde_{3}(\alpha) \vert^2 \frac{\vert \Etilde_{3}(\alpha) \vert}{\vert z \vert^{1/3}}
+
\frac{\vert \Etilde_{3}(\alpha) \vert}{\vert z \vert}
+
N^{2/3}\vert \Etilde_{3}(\alpha) \vert^2.
\end{align}
 Using  \eqref{S-tilde-uncond-approx} and \eqref{z-estim} we get
 \begin{align}
 \notag
J_2
& \ll
H 
\int_{-B/H}^{B/H}
\vert \Stilde_{3}(\alpha) \vert^2 \frac{\vert \Etilde_{3}(\alpha) \vert}{\vert z \vert^{1/3}}  \, \dx \alpha
+ 
H 
\int_{-B/H}^{B/H}
\frac{\vert \Etilde_{3}(\alpha) \vert}{\vert z \vert}  \, \dx \alpha
+
HN^{2/3}
\int_{-B/H}^{B/H}
\vert \Etilde_{3}(\alpha) \vert^2  \, \dx \alpha
\\
\label{J2-uncond-split}&
= H(  E_1+  E_2  +N^{2/3}E_3 ),
 \end{align}
 say. 
 By  \eqref{UH-estim} and  Lemma \ref{LP-Lemma-gen} we obtain
 that,  for every $\eps>0$, there exists $c_1=c_1(\eps)>0$ such that
 \begin{equation}
 \label{E3-estim-series}
E_3 
\ll  
   N^{-1/3} \exp \Big( - c_{1}  \Big( \frac{L}{\log L} \Big)^{1/3} \Big)  
\end{equation}
provided that $B/H \le N^{-13/18 - \eps}$, \emph{i.e.}, 
$H \ge BN^{ 13/18+ \eps}$.
 By the Cauchy-Schwarz inequality, \eqref{z-estim} 
 and \eqref{E3-estim-series} we obtain
  that,  for every $\eps>0$, there exists $c_1=c_1(\eps)>0$ such that
 \begin{equation}
\label{E2-estim-series}
E_2
  \ll 
 E_3^{1/2}   \Bigl(\int_{-B/H}^{B/H}
\frac{\dx \alpha}{\vert z \vert^2}  
\Bigr)^{1/2}
 \ll
 E_3^{1/2}  N^{1/2}
 \ll
    N^{1/3} \exp \Big( - \frac{c_{1}}{2}  \Big( \frac{L}{\log L} \Big)^{1/3} \Big), 
\end{equation}
provided that  $H \ge BN^{ 13/18+ \eps}$.
 By using twice the Cauchy-Schwarz inequality, Lemma \ref{Hua-lemma-series},
  \eqref{z-estim} 
 and \eqref{E3-estim-series} we obtain
  that,  for every $\eps>0$, there exists $c_1=c_1(\eps)>0$ such that
 \begin{align}
\notag
E_1
&  \ll 
 E_3^{1/2}   \Bigl(\int_{-B/H}^{B/H}
\frac{\vert \Stilde_{3}(\alpha) \vert^4}{\vert z \vert^{2/3}}  \ \dx \alpha
\Bigr)^{1/2}
\ll
 E_3^{1/2}   
\Bigl(\int_{-1/2}^{1/2}
\vert \Stilde_{3}(\alpha) \vert^8   \ \dx \alpha
\Bigr)^{1/4}
\Bigl(\int_{-B/H}^{B/H}
\frac{\dx \alpha}{\vert z \vert^{4/3}}  
\Bigr)^{1/4}
\\&
\label{E1-estim-series}
 \ll
  E_3^{1/2}  N^{1/2}L^{A/4}
  \ll
    N^{1/3} \exp \Big( - \frac{c_{1}}{4}  \Big( \frac{L}{\log L} \Big)^{1/3} \Big), 
\end{align}
provided that  $H \ge BN^{ 13/18+ \eps}$.
Hence by \eqref{J2-uncond-split}-\eqref{E1-estim-series}
we finally can write   that,  for every $\eps>0$, there exists $c_1=c_1(\eps)>0$ such that
 \begin{equation}
\label{J2-estim-series}
J_2
\ll
H N^{1/3} \exp \Big( - \frac{c_{1}}{4} \Big( \frac{L}{\log L} \Big)^{1/3} \Big) ,
 \end{equation}
 provided that
 $H \ge BN^{ 13/18+ \eps}$.
Summing up, by \eqref{I1-split}-\eqref{J1-eval-series} and \eqref{J2-estim-series} we have   that,  for every $\eps>0$, there exists $c_1=c_1(\eps)>0$ such that
\begin{equation}
  \label{I1-final-eval}
I_1 =
\Gamma\Bigl(\frac{4}{3}\Bigr)^3
\frac{HN^{1/3}}{e}
+\Odig{ 
H N^{1/3} \exp \Big( - \frac{c_{1}}{4} \Big( \frac{L}{\log L} \Big)^{1/3} \Big) 
}
\end{equation}
 provided that
 $H \ge BN^{ 13/18+ \eps}$.
\subsection{Final words}

Summing up,  
by    \eqref{main-dissection-series}-\eqref{I2-estim}, 
\eqref{I3-estim-series} and \eqref{I1-final-eval}
we have  that,  for every $\eps>0$, there exists $c_1=c_1(\eps)>0$ such that
\begin{align*} 
\sum_{n=N+1}^{N+H} 
e^{-n/N}
r(n)
&=
\Gamma\Bigl(\frac{4}{3}\Bigr)^3
\frac{HN^{1/3}}{e}
+\Odig{H N^{1/3} \exp \Big( - \frac{c_{1}}{4} \Big( \frac{L}{\log L} \Big)^{1/3} \Big) 
+
\frac{H }{B }  N^{1/3+\eps}}
\end{align*}
provided that  $H \ge BN^{ 13/18+ \eps}$.
The second error term is dominated by the first one 
since $B=N^{2\eps}$ by   \eqref{B-def}. Hence
we can write  that,  for every $\eps>0$, there exists $C=C(\eps)>0$ such that
\begin{equation}
\label{almost-done-2}
\sum_{n=N+1}^{N+H} 
e^{-n/N}
r(n)
=
\Gamma\Bigl(\frac{4}{3}\Bigr)^3
\frac{HN^{1/3}}{e}
+\Odig{H N^{1/3} \exp \Big( - C\Big( \frac{L}{\log L} \Big)^{1/3} \Big)}
\end{equation}
provided that $H \ge N^{13/18+ 3\eps}$.
{}From  $e^{-n/N}=e^{-1}+ \Odi{H/N}$ for $n\in[N+1,N+H]$, $1\le H \le N$,
we get  that,  for every $\eps>0$, there exists $C=C(\eps)>0$ such that
\begin{align*}
   \sum_{n = N+1}^{N + H} 
r(n)
&=
\Gamma\Bigl(\frac{4}{3}\Bigr)^3
HN^{1/3} 
+ \Odig{H N^{1/3} \exp \Big( - C\Big( \frac{L}{\log L} \Big)^{1/3} \Big)}
+
  \Odig{\frac{H}{N}\sum_{n = N+1}^{N + H} r(n)
}
\end{align*}
provided that $H \ge N^{13/18+ 3\eps}$ and $H \le N$.
Using $e^{n/N}\le  e^{2}$ 
and \eqref{almost-done-2},
 the last error term is
$\ll H^2N^{-2/3}$.
Hence we get  that,  for every $\eps>0$, there exists $C=C(\eps)>0$ such that
\begin{equation*}
   \sum_{n = N+1}^{N + H} 
r(n)
=
\Gamma\Bigl(\frac{4}{3}\Bigr)^3
HN^{1/3} 
+\Odig{H N^{1/3} \exp \Big( - C \Big( \frac{L}{\log L} \Big)^{1/3} \Big)},
\end{equation*}
provided that $N^{13/18+ 3\eps} \le  H \le N^{1-\eps}$.
Theorem \ref{thm-uncond} follows by rescaling $\eps$.

\section{The conditional case}
\label{conditional}
From now on we assume the Riemann Hypothesis holds. 
Comparing with section \ref{unconditional} we can simplify the setting.
Recalling \eqref{r-def} and  $\Gamma(4/3)= (1/3)\Gamma(1/3)$,
we have  
\begin{align}
 \notag
\sum_{n=N+1}^{N+H} e^{-n/N}  r(n)
& =  
\int_{-1/2}^{1/2}\Vtilde_{3}(\alpha)^4    U(-\alpha,H)e(-N\alpha) \, \dx \alpha  
  \\
 \notag
 &
 = \Gamma \Bigl(\frac{4}{3}\Bigr)^4
 \int_{-1/2}^{1/2}\frac{U(-\alpha,H)}{z^{4/3}}e(-N\alpha) \, \dx \alpha 
 \\
 \notag
 &
\hskip1cm
+
\int_{-1/2}^{1/2} \Bigl( \Stilde_{3}(\alpha)^4  - \frac{\Gamma(4/3)^4}{z^{4/3}}\Bigr)  
 U(-\alpha,H)e(-N\alpha) \, \dx \alpha 
 \\
 \notag
 &
\hskip1cm
+
\int_{-1/2}^{1/2} \bigl( \Vtilde_{3}(\alpha)^4  -\Stilde_{3}(\alpha)^4 \bigr)  
 U(-\alpha,H)e(-N\alpha) \, \dx \alpha 
\\
 \label{main-dissection-RH-series}
 & 
= \J_1 +\J_2 + \J_3,
 \end{align}
 say.
 Now we evaluate these terms.
 
\subsection{Evaluation of $\J_1$}

By Lemma \ref{Laplace-formula},
a direct calculation  gives
\begin{align}
\notag
\J_1
&
=
\Gamma\Bigl(\frac{4}{3}\Bigr)^3\
\sum_{n=N+1}^{N+H} e^{-n/N} n^{1/3}
+\Odig{\frac{H}{N}}
=
\frac{\Gamma(4/3)^3}{e}\
\sum_{n=N+1}^{N+H} n^{1/3}
+\Odig{\frac{H}{N}+\frac{H^2}{N^{2/3}}}
\\
\label{J1-eval-RH-series}
&
=
\Gamma\Bigl(\frac{4}{3}\Bigr)^3
\frac{HN^{1/3}}{e}
+\Odig{\frac{H^2}{N^{2/3}}+N^{1/3}}.
 \end{align}

\subsection{Estimate of $\J_2$}
Recall that
\( 
\Etilde_{3}(\alpha) : =\Stilde_3(\alpha) - \frac{\Gamma(4/3)}{z^{1/3}}.
\)
Using $  f^2 - g^2 = 2g (f-g) + (f-g)^2$
  we can write that
\begin{align*}
\Stilde_{3}(\alpha)^4  - \frac{\Gamma(4/3)^4}{z^{4/3}} 
&=
\Bigl( 
\Stilde_{3}(\alpha) ^2  +  
\frac{\Gamma(4/3)^2}{z^{2/3}}
\Bigr)
\Bigl(
\Stilde_{3}(\alpha) ^2  - 
\frac{\Gamma(4/3)^2}{z^{2/3}}
\Bigr)
\\&
\ll
(\vert \Stilde_{3}(\alpha) \vert^2 + \vert z \vert^{-2/3})
(
 \vert z \vert^{-1/3}   \vert \Etilde_{3}(\alpha) \vert
+
 \vert \Etilde_{3}(\alpha) \vert^2
 ).
\end{align*}
Hence
\begin{align}
\notag
\J_2 
&\ll
\int_{-1/2}^{1/2}
\frac{\vert \Stilde_{3}(\alpha) \vert^2}{\vert z \vert^{1/3}} \vert \Etilde_{3}(\alpha) \vert
\vert U(-\alpha,H)\vert \, \dx \alpha
+
\int_{-1/2}^{1/2} 
\frac{\vert \Etilde_{3}(\alpha) \vert}{\vert z \vert}\vert U(-\alpha,H)\vert \, \dx \alpha
\\ \notag &
\hskip 1cm+
\int_{-1/2}^{1/2}
(\vert \Stilde_{3}(\alpha) \vert^2 + \vert z \vert^{-2/3})
 \vert \Etilde_{3}(\alpha) \vert^2 \vert U(-\alpha,H)\vert \, \dx \alpha
\\
\label{J2-RH-split}
&
 =
\I_1+\I_2+\I_3, 
\end{align}
say. Let
\[
\E :=  \int_{-1/2}^{1/2}
\vert \Etilde_{3}(\alpha) \vert^2 \vert U(-\alpha,H) \vert \, \dx \alpha.
\]
By \eqref{UH-estim}, Lemma \ref{LP-Lemma-gen}
and a partial integration we obtain
\begin{align}
\notag
\E 
&\ll
H  
\int_{-1/H}^{1/H}
\vert  \Etilde_{3}(\alpha) \vert^2   \, \dx \alpha 
+
\int_{1/H}^{1/2}
\frac{\vert  \Etilde_{3}(\alpha) \vert^2}{\alpha} \, \dx \alpha 
\\& \label{E-estim}
\ll 
H  \frac{N^{1/3}  L^2}{H} 
+   
N^{1/3} L^2
+
\int_{1/H}^{1/2}  
\frac{ N^{1/3}  L^2} {\xi} \ \dx \xi 
\ll N^{1/3} L^3.
\end{align}
By \eqref{z-estim},  $\Stilde_{3}(\alpha)  \ll  N^{1/3}$
and  \eqref{E-estim} we obtain
 \begin{equation}
 \label{I3-estim-RH-series}
\I_3 
\ll  
N^{2/3} \E 
\ll 
N L^3.
\end{equation}
 By the Cauchy-Schwarz inequality,  \eqref{z-estim}-\eqref{UH-estim}
 and   \eqref{E-estim},  we obtain
 \begin{align}
 \notag
\I_2
&
\ll 
\E^{1/2}
\Bigl(
\int_{-1/2}^{1/2}
 \frac{\vert U(-\alpha,H) \vert }{\vert z \vert^{2}} \, \dx \alpha
\Bigr)^{1/2}
\ll
\E^{1/2}
\Bigl(
 H N^{2}
\int_{-1/N}^{1/N}  \dx \alpha
 +
 H
\int_{1/N}^{1/H}
 \frac{ \dx \alpha}{  \alpha  ^{2}}   
 +
 \int_{1/H}^{1/2}
 \frac{ \dx \alpha}{  \alpha  ^{3}}  
\Bigr)^{1/2}
\\
\label{I2-estim-RH-series}
& \ll
H^{1/2} N^{2/3} L^{3/2}.
\end{align}
By the Cauchy-Schwarz inequality we obtain
 \[
 \I_1
\ll 
\E^{1/2}
\Bigl(
\int_{-1/2}^{1/2}
\vert \Stilde_3(\alpha) \vert^4 \frac{\vert U(-\alpha,H) \vert }{\vert z \vert^{2/3}} \, \dx \alpha
\Bigr)^{1/2}.
 \]
Again by the Cauchy-Schwarz inequality,  \eqref{z-estim}-\eqref{UH-estim}
 and \eqref{E-estim},  we obtain
 \begin{align}
 \notag
\I_1
&
\ll 
\E^{1/2}
\Bigl(
\int_{-1/2}^{1/2}
\vert \Stilde_3(\alpha) \vert^8  \, \dx \alpha
\Bigr)^{1/4}
\Bigl(
\int_{-1/2}^{1/2}
\frac{\vert U(-\alpha,H) \vert^2 }{\vert z \vert^{4/3}} \, \dx \alpha
\Bigr)^{1/4}
\\
\notag
&
\ll
N^{1/6} L^{3/2} N^{5/12} L^{A/4}
\Bigl(
 H^2 N^{4/3}
\int_{-1/N}^{1/N}  \dx \alpha
 +
 H^2
\int_{1/N}^{1/H}
 \frac{ \dx \alpha}{\alpha^{4/3}}   
 +
 \int_{1/H}^{1/2}
 \frac{ \dx \alpha}{\alpha^{10/3}}  
\Bigr)^{1/4}
\\
\label{I1-estim-RH-series}
& \ll
H^{1/2} N^{2/3} L^{3/2+A/4}.
\end{align}
Summing up by \eqref{J2-RH-split} and \eqref{I3-estim-RH-series}-\eqref{I1-estim-RH-series},
we  can finally write that
 \begin{equation}
\label{J2-estim-RH-series}
\J_2 \ll   H^{1/2} N^{2/3} L^{3/2+A/4} + N L^3.
 \end{equation}

\subsection{Estimate of $\J_3$} %

It is clear that $\J_3=I_3$ of section \ref{I3-estim}. Hence by \eqref{I3-estim-series}
we obtain
\begin{equation}
 \label{J3-estim-RH-series}  
\J_3  
\ll 
H^{3/4}N^{5/12+\eps}.
\end{equation}

\subsection{Final words}

Summing up,  
by    \eqref{main-dissection-RH-series}-\eqref{J1-eval-RH-series}, 
\eqref{J2-estim-RH-series} and \eqref{J3-estim-RH-series},
there exists $B=B(A)>0$ such that
we have
\begin{equation}
\label{almost-done}
\sum_{n=N+1}^{N+H} 
e^{-n/N}
r(n)
=
\Gamma\Bigl(\frac{4}{3}\Bigr)^3
\frac{HN^{1/3}}{e}
+\Odig{\frac{H^2}{N^{2/3}}
+ H^{3/4}N^{5/12+\eps}
+ H^{1/2} N^{2/3} L^{B}
+ NL^3}
\end{equation}
which is an asymptotic formula 
$\infty(N^{2/3}L^{2B}) \le H \le \odi{N}$.
{}From  $e^{-n/N}=e^{-1}+ \Odi{H/N}$ for $n\in[N+1,N+H]$, $1\le H \le N$,
we get 
\begin{align*}
   \sum_{n = N+1}^{N + H} 
r(n)
&=
\Gamma\Bigl(\frac{4}{3}\Bigr)^3
HN^{1/3} 
+\Odig{\frac{H^2}{N^{2/3}}
+ H^{3/4}N^{5/12+\eps}
+ H^{1/2} N^{2/3} L^{B}
+NL^3} 
+
  \Odig{\frac{H}{N}\sum_{n = N+1}^{N + H} r(n)
}.
\end{align*}
Using $e^{n/N}\le  e^{2}$ 
and \eqref{almost-done},
 the last error term is
$\ll H^2N^{-2/3}+  H^{7/4}N^{-7/12+\eps} + H^{3/2} N^{-2/3} L^{B} + HL^3$.
Hence we get
\begin{equation*}
   \sum_{n = N+1}^{N + H} 
r(n)
=
\Gamma\Bigl(\frac{4}{3}\Bigr)^3
HN^{1/3} 
+\Odig{\frac{H^2}{N^{2/3}}+ H^{3/4}N^{5/12+\eps}
+ H^{1/2} N^{2/3} L^{B}
+NL^3} ,
\end{equation*}
uniformly for $\infty(N^{2/3}L^{2B}) \le H \le \odi{N}$, $B>3/2$.
Theorem \ref{thm-RH} follows.
 
\renewcommand{\bibliofont}{\normalsize}
  
%
 
\providecommand{\bysame}{\leavevmode\hbox to3em{\hrulefill}\thinspace}
\providecommand{\MR}{\relax\ifhmode\unskip\space\fi MR }
\providecommand{\MRhref}[2]{%
  \href{http://www.ams.org/mathscinet-getitem?mr=#1}{#2}
}
\providecommand{\href}[2]{#2}

\vskip0.5cm
\noindent
\begin{tabular}{l@{\hskip 20mm}l}
Alessandro Languasco               & Alessandro Zaccagnini\\
Universit\`a di Padova     & Universit\`a di Parma\\
 Dipartimento di Matematica  &  Dipartimento di Scienze Matematiche,\\
 ``Tullio Levi-Civita'' &  Fisiche e Informatiche \\
Via Trieste 63                & Parco Area delle Scienze, 53/a \\
35121 Padova, Italy            & 43124 Parma, Italy\\
{\it e-mail}: alessandro.languasco@unipd.it      & {\it e-mail}:
alessandro.zaccagnini@unipr.it  
\end{tabular}
 
 \end{document}